\documentclass[10 pt]{amsart}
\usepackage{amsmath}
\usepackage{amssymb}
\usepackage{amsfonts}

\theoremstyle{plain}
\newtheorem{theorem}  {Theorem} [section]
\newtheorem{lemma}    [theorem] {Lemma} 
\newtheorem{corollary}[theorem] {Corollary}

\theoremstyle{definition}
\newtheorem{definition}[theorem]{Definition}

\theoremstyle{remark}
\newtheorem{remark}    [theorem]{Remark}
\newtheorem{example}   [theorem]{\bf Example}

\begin{document}

\title
{$\mathbf{{L_{\infty}}}$ structures on spaces with 3 one-dimensional components}
\author{Marilyn Daily}
\address{Department of Mathematics, North Carolina State University,
Raleigh NC 27695}
\email{medaily@eos.ncsu.edu}
\maketitle

\section{Introduction}

$L_{\infty}$ structures have been a subject of recent interest in physics, where they occur in 
closed string theory \cite{Z93} and in gauge theory \cite{FLS01}.
This paper provides a class of easily constructible examples of $L_n$ and $L_{\infty}$ structures
on graded vector spaces with 3 one-dimensional components.
In particular, we will show how to classify {\bf all} possible $L_n$ and 
$L_{\infty}$ structures on 
$V = V_m \oplus V_{m+1} \oplus V_{m+2}$
when $dim(V_m)=dim(V_{m+1})=dim(V_{m+2})=1$.  
We will also give necessary and sufficient conditions under which a space with an $L_3$ structure
is a differential graded Lie algebra.
In the process, we will show that some of these d.g. Lie algebras possess a 
nontrivial $L_n$ structure for higher $n$.

The fundamental properties of $L_{\infty}$ structures (including many of the definitions below)
were developed earlier in \cite{LM95} and \cite{LS93}.
The results in this paper are a portion of my PhD research 
under the supervision of Tom Lada.
I would also like to thank Jim Stasheff for many helpful suggestions and comments.

\section{Basic definitions and notation}

Throughout this paper, it is assumed that $V = V_m \oplus V_{m+1} \oplus V_{m+2}$
is a graded vector space over a fixed field of characteristic zero.
Since we are working with graded vector spaces, the Koszul sign
convention will be employed: whenever two objects of
degrees $p$ and $q$ are commuted, a factor of $(-1)^{pq}$ is
introduced.  For a permutation $\sigma$ acting on a string of symbols,
let $\varepsilon(\sigma)$ denote the total effect of these signs; then
$\chi(\sigma) = (-1)^{\sigma}\varepsilon(\sigma)$ where $(-1)^{\sigma}$
is the sign of the permutation $\sigma$.

\begin {remark}
$\chi(\sigma)$ has the following properties, which can be easily verified:
\begin{enumerate}
\item
If a permutation $\sigma$ just reorders a string of consecutive elements of odd degree, then $\chi(\sigma)=1$.
\item
If a permutation $\sigma$ just reorders a string of consecutive elements of even degree, then $\chi(\sigma)=(-1)^{\sigma}$.
\item
If a permutation $\sigma$ just transposes two (arbitrarily positioned) elements of even degree, then
$\chi(\sigma)=-1$.
\item
If a permutation $\sigma$ moves an element of even degree past a string of $k$ elements of arbitrary degree, then $\chi(\sigma)=(-1)^k$.
\item
If a permutation $\sigma$ moves an element of odd degree past a string of $k$ elements of even degree, then $\chi(\sigma)=(-1)^k$.
\item
If a permutation $\sigma$ moves a string of $j$ even elements past a string of $k$ elements of arbitrary degree, then $\chi(\sigma)=(-1)^{jk}$.
\end{enumerate}
\end {remark}

\begin {definition} A $(p,n-p)$ {\em unshuffle} is a permutation $\sigma\in S_n$ such that
$$\sigma (1)< \sigma (2)<\dots <\sigma(p) \text{ and } \sigma (p+1) < \sigma (p+2) <\dots <\sigma(n).$$
\end {definition}

\begin {definition} A map $l_n:V^{\otimes n}\rightarrow V$ has {\em degree} $k$ if
$l_n(x_1\otimes x_2\otimes\dots\otimes x_n) \subset V_N$, where $N=k+\sum_{i=1}^n \vert x_i \vert.$
\end {definition}

\begin {definition} A linear map $l_n:V^{\otimes n}\rightarrow V$ is {\em skew symmetric} if
$$
l_n(x_{\pi(1)}\otimes\dots\otimes x_{\pi(n)})=\chi(\pi)l_n(x_1\otimes\dots\otimes x_n)
\ \forall\ \pi\in S_n.
$$
\end {definition}

\begin {definition} A skew symmetric linear map $l:V^{\otimes n}\rightarrow V$ can be extended to a map
$l:V^{\otimes n+k}\rightarrow V^{\otimes 1+k}$ by the rule,
$$
l(x_{1}\otimes\dots\otimes x_{n+k})=
\sum_{\substack{\sigma \text{ is } (n,k) \\ \text{ unshuffle}}}
\chi(\sigma)l(x_{\sigma(1)}\otimes\dots\otimes x_{\sigma(n)})\otimes x_{\sigma(n+1)}\otimes\dots\otimes x_{\sigma(n+k)}.
$$
\end {definition}

\begin {remark}
The extension of a skew linear map need not be skew.  For example:  
$l_1(v\otimes u\otimes w)\neq -(-1)^{uv}l_1(u\otimes v\otimes w)$.
However, it can be shown that the composition
$l_{n-p+1}\circ l_p:V^{\otimes n}\rightarrow V$ is a skew symmetric linear operator.
\end {remark}

\begin {definition} Let $\mathcal J_n:V^{\otimes n}\rightarrow V$ be defined by
$$
\mathcal J_n (x_1\otimes x_2\otimes\dots\otimes x_n)
=\sum_{p=1}^n
(-1)^{p(n-p)}l_{n-p+1}\circ l_p (x_1\otimes x_2\otimes\dots\otimes x_n).
$$
(using the notation of extended maps defined above).  
\end {definition}

\begin {remark}For ease of computation,
$\mathcal J_n(x_1\otimes x_2\otimes\dots\otimes x_n)$ can be rewritten as
$$
\sum_{p=1}^n (-1)^{p(n-p)}
\sum_{\substack{ \sigma \text{ is } (p,n-p) \\ \text{unshuffle} }}\chi(\sigma)
l_{n-p+1}(l_p(x_{\sigma(1)}\otimes\dots\otimes x_{\sigma(p)})\otimes x_{\sigma(p+1)}\otimes\dots\otimes x_{\sigma(n)}).
$$
\end {remark}

\begin {definition} The {\em nth generalized Jacobi identity} is the equation $$\mathcal J_n = 0.$$
\end {definition}

\begin {remark}
If we denote $l_2(x\otimes y)=[x,y]$ and consider the special case when $n=3$ and $l_3=0$, the generalized Jacobi identity becomes
the more familiar graded Jacobi identity, since
$$
\sum_{\substack{ \sigma \text{ is } (2,1) \\ \text{ unshuffle} }}
\chi(\sigma )
l_2(l_2(x_{\sigma(1)}\otimes x_{\sigma(2)})\otimes x_{\sigma(3)}) =0 
$$
can be rewritten in bracket notation as
$$
\sum_{\substack{ \sigma \text{ is } (2,1) \\ \text{ unshuffle} }}
\chi(\sigma )
[ [ x_{\sigma(1)}, x_{\sigma(2)})], x_{\sigma(3)}] =0 
$$
which then expands to the usual graded Jacobi identity,
$$
[ [x_1,x_2],x_3] - (-1)^{\vert x_2 \vert \vert x_3 \vert} [ [x_1,x_3],x_2]
+ (-1)^{\vert x_1\vert (\vert x_2 \vert + \vert x_3 \vert)}[ [x_2,x_3], x_1] = 0.
$$
\end {remark}

\begin {remark}
Note that $l_{n-p+1}\circ l_p:V^{\otimes n}\rightarrow V$ is a skew linear operator of degree $3-n$.
Therefore, $\mathcal J_n$ is also a skew linear operator of degree $3-n$.
In particular, whenever an element of even degree is repeated, the Jacobi identity is automatically satisfied. 
\end {remark}

\begin {definition} An $L_m$ {\em structure} is a graded vector space
$V$ endowed with a collection of linear maps $\{l_k:V^{\otimes k}\rightarrow V,\  1 \leq k \leq m\}$ 
with $deg(l_k)=2-k$ which are skew symmetric and satisfy each Jacobi identity $\mathcal J_n = 0$, $1\leq n\leq m$.
\end {definition}

\begin {remark}
This is the \lq\lq cochain version\rq\rq\ of the definition of an $L_m$ structure,  having $deg(l_k)=2-k$.
The definition for chain complexes is the same, but with $deg(l_k)=k-2$.
We will be working exclusively with the cochain version in this paper.
\end {remark}

\begin {definition} An $L_{\infty}$ {\em structure} is a graded vector space
$V$ endowed with a collection of linear maps $\{l_k:V^{\otimes k}\rightarrow V,\ k \in \mathbb N \}$ 
such that $\mathcal J_n = 0 \ \forall\ n \in \mathbb N$.
\end {definition}

\begin {definition}
A {\em differential graded (d.g.) Lie algebra} is an $L_2$ structure which satisfies the identity,
$$\sum_{\substack{ \sigma \text{ is } (2,1) \\ \text{ unshuffle} }}
\chi(\sigma )
l_2(l_2(x_{\sigma(1)}\otimes x_{\sigma(2)})\otimes x_{\sigma(3)}) =0.$$
\end {definition}

%
%

\section {$\mathbf{V=V_0\oplus V_1\oplus V_2}$ where each is 1-dimensional}

This section includes necessary and sufficient conditions 
under which skew linear operators form an $L_m$ or $L_{\infty}$ structure on 
$V=V_0\oplus V_1\oplus V_2$,
as well as necessary and sufficient conditions for an $L_3$ structure
on this space to be a differential graded Lie algebra.

Specifically,  consider $V=\oplus V_i$ where
\begin {alignat*} {4}
V_0 &= <v> \qquad & \qquad V_1 &= <w> \qquad & \qquad V_2 &= <x> \qquad & \qquad V_k &= {0},\ k \notin \{ 0,1,2 \}
\end {alignat*}
Define skew linear operators $l_n:V^{\otimes n}\rightarrow V$ of degree $2-n$ by
\begin {equation}
\left\{\ 
\begin {aligned}
l_n(v\otimes w^{\otimes n-1})&=a_n w          \!\quad \forall\ n\geq 1 \\
l_n(w^{\otimes n})&=b_n x                      \,\quad \forall\ n\geq 1 \\
l_n(v\otimes w^{\otimes n-2}\otimes x)&=c_n x  \,\quad \forall\ n\geq 2 \\
l_n(w^{\otimes n-1}\otimes x)&=0              \,\qquad \forall\ n       \\
l_n(v^{\otimes i}\otimes w^{\otimes n-i-j}\otimes x^{\otimes j})&=0 \!\qquad \ \forall\ i,\ j>1
\end {aligned}
\right.
\label {Tag1}
\end {equation}
where $a_n$, $b_n$, $c_n$ are constants.
Since $c_n$ is defined above only for $n\geq 2$, we can define $c_1=0$ for the sake of convenience.

The goal in this section is to show that this space
can have just two types of $L_{\infty}$ structure:
\begin {itemize}
\item
$b_p=0\ \forall\ p\in \mathbb N$ but all constants $a_n$, $c_n$ are arbitrarily chosen
\item
If some $b_p$ is nonzero, the remaining constants $b_n$ and all constants $c_n$
can be arbitrarily chosen, but all constants $a_n$ are then uniquely determined by those
choices.
\end {itemize}

\begin {remark}The above list 
\eqref{Tag1}
includes all possible skew linear operators of degree $2-n$, by the 
following argument:
\begin {itemize}
\item 
$l_n(v^{\otimes i}\otimes w^{\otimes n-i-j}\otimes x^{\otimes j})=0 \quad \forall\ i,\ j >1$
because a repeated component of even degree makes a skew operator zero.
\item
$l_n(w^{\otimes n-1}\otimes x)\subset V_{n-1+2+2-n}=V_{3} $ which is zero.
\end {itemize}
The action of $l_k$ on the remaining generators is determined by its skewness.
\end{remark}

In order to determine when these operators define an $L_m$ structure on $V$,
we need to know when $\mathcal J_n = 0 \ \forall\ n \leq m$.  The following two lemmas establish
the conditions which force $\mathcal J_n = 0$:

\begin {lemma} \label {lem3.2}
Given $V=V_0\oplus V_1\oplus V_2$ with the skew linear operators defined in \eqref{Tag1},
$\mathcal J_n = 0$ if and only if 
$\mathcal J_n(v\otimes w^{\otimes n-1}) = 0$. 
\end {lemma} 

\begin {proof}
Since $deg(\mathcal J_n) = 3-n$, 
$$
\mathcal J_n(V_0^{\otimes i}\otimes V_1^{\otimes n-i-j}\otimes V_2^{\otimes j}) 
\subset V_{n-i-j+2j + 3 - n} = V_{-i+j+3}.
$$
Thus
$\mathcal J_n(V_0^{\otimes i}\otimes V_1^{\otimes n-i-j}\otimes V_2^{\otimes j})$ 
can only be nonzero if $j-i+3\in \{ 0,1,2 \}$
(since these are the only nonzero components of $V$).
Equivalently,
$\mathcal J_n(V_0^{\otimes i}\otimes V_1^{\otimes n-i-j}\otimes V_2^{\otimes j})$ 
can only be nonzero if $j-i\in \{ -3,-2,-1 \}$.  
But 
$\mathcal J_n(V_0^{\otimes i}\otimes V_1^{\otimes n-i-j}\otimes V_2^{\otimes j})=0$ 
whenever $i>1$ (since $\mathcal J_n$ is skew and a parameter of even degree would be repeated).
Thus 
$\mathcal J_n(V_0^{\otimes i}\otimes V_1^{\otimes n-i-j}\otimes V_2^{\otimes j})$ 
can only be nonzero when $i = 1$ and $j = 0$.
\end {proof}


\begin {lemma} \label {lem3.3}
Given $V=V_0\oplus V_1\oplus V_2$ with the operators defined in \eqref{Tag1},
$$
\mathcal J_n(v\otimes w^{\otimes n-1})
 = 
\sum_{p=1}^n (-1)^{p(n-p)}(-1)^n b_p \left[ (-1) \tbinom{n-1}{n-p} a_{n-p+1} 
+ \tbinom{n-1}{p} c_{n-p+1} \right] x.
$$
\end {lemma}

\begin {proof}
Before calculating $l_{n-p+1}\circ l_p(v\otimes w^{\otimes n-1})$, it is helpful to consider the
ways in which a $(p,n-p)$ unshuffle of $v\otimes w^{\otimes n-1}$ could possibly rearrange the terms.
Here, there are only two possibilities, since ordering guarantees that the $v$ must
be placed in either position $1$ or position $p+1$.

There are $\tbinom{n-1}{p-1}$ unshuffles which place $v$ in position $1$.  
In each case, $\chi(\sigma)=+1$,
since each of these permutations can be accomplished solely by commuting elements of odd degree.

Similarly, there are $\tbinom{n-1}{p}$ unshuffles which place $v$ in position $p+1$.  
In each case, $\chi(\sigma)=(-1)^p$.
This is because we can reorder the consecutive $w$ terms (which does not change the sign), and
then simply shift $v$ past $p$ $w$ terms.

\begin {align*}
l_{n-p+1}\circ l_p(v\otimes w^{\otimes n-1})
&=
\tbinom{n-1}{p-1}l_{n-p+1}(l_p(v\otimes w^{\otimes p-1})\otimes w^{\otimes n-p})
\\
&\quad +
\tbinom{n-1}{p}(-1)^{p}l_{n-p+1}(l_p(w^{\otimes p})\otimes v \otimes w^{\otimes n-p-1})
\\
&=
\tbinom{n-1}{p-1}l_{n-p+1}(a_p w\otimes w^{\otimes n-p})
\\
&\quad +
\tbinom{n-1}{p}(-1)^{p}(-1)^{n-p}l_{n-p+1}(b_p v\otimes w^{\otimes n-p-1}\otimes x)
\\
&=
\tbinom{n-1}{p-1} a_p b_{n-p+1}x
 +
(-1)^{n}\tbinom{n-1}{p} b_p c_{n-p+1}x.
\end {align*}

Since $\mathcal J_n = 
\sum_{p=1}^n (-1)^{p(n-p)} l_{n-p+1}\circ l_p$, we have

$$
\mathcal J_n(v\otimes w^{\otimes n-1})
= 
\sum_{p=1}^n (-1)^{p(n-p)} \left[ \tbinom{n-1}{p-1} a_p b_{n-p+1}x
+ (-1)^n \tbinom{n-1}{p} b_p c_{n-p+1}x \right].
$$

Thus $\mathcal J_n(v\otimes w^{\otimes n-1})$ is equal to 
$$
\sum_{p=1}^n (-1)^{p(n-p)} \tbinom{n-1}{p-1} a_p b_{n-p+1}x
+ 
\sum_{p=1}^n (-1)^{p(n-p)}
(-1)^n \tbinom{n-1}{p} b_p c_{n-p+1}x.
$$
Substituting $q=n-p+1$ in the first sum and re-indexing, this is equivalent to
$$
\sum_{q=1}^n (-1)^{(n-q+1)(q-1)} \tbinom{n-1}{n-q} a_{n-q+1} b_{q} x
+ \sum_{p=1}^n (-1)^{p(n-p)}(-1)^n \tbinom{n-1}{p} b_p c_{n-p+1} x.
$$
Recombining the two sums, this is equivalent to 
$$
\sum_{p=1}^n (-1)^{p(n-p)}(-1)^n b_p \left[ (-1) \tbinom{n-1}{n-p} a_{n-p+1} 
+ \tbinom{n-1}{p} c_{n-p+1} \right] x.
$$
\end {proof}

\begin {remark}
Thus in $V=V_0\oplus V_1\oplus V_2$, 
$$
\mathcal J_n = 0 \Longleftrightarrow 0=
\sum_{p=1}^n (-1)^{p(n-p)}(-1)^n b_p \left[ (-1) \tbinom{n-1}{n-p} a_{n-p+1} 
+ \tbinom{n-1}{p} c_{n-p+1} \right].
$$
In particular, $v$ is $L_1$ (i.e. a differential graded vector space) if and only if $0=a_1 b_1$.
\end {remark}


\begin {lemma} \label {lem3.4}
Suppose $V=V_0\oplus V_1\oplus V_2$ with the skew linear operators defined in \eqref{Tag1}
has an $L_{\infty}$ structure and $a_1 \neq 0$.  Then $b_n = 0 \ \forall\  n \in \mathbb N$. 
\end {lemma}

\begin {proof} 
Suppose $a_1 \neq 0$.  Then $b_1 = 0$ (since the first Jacobi identity requires that $0=a_1 b_1$).
Now we do an induction proof:  Suppose $b_k=0 \ \forall\ k < n$.  
The $n$th Jacobi identity requires that
$$
0=\sum_{p=1}^n (-1)^{p(n-p)}(-1)^n b_p \left[ (-1) \tbinom{n-1}{n-p} a_{n-p+1} + \tbinom{n-1}{p} c_{n-p+1} \right].
$$
Since we're supposing that $b_k = 0\ \forall\ k < n$, this collapses to 
$0=(-1)^n b_n (n-1) a_{1}$.
But we also supposed that $a_1 \neq 0$, which forces $b_n = 0$.
\end {proof}


If some constant $b_k$ is nonzero, then the following lemma proves that the constraint
$\mathcal J_n=0$ is equivalent to an explicit formula for $a_{n-k+1}$:

\begin {lemma} \label {lem3.5}
Given $V=V_0\oplus V_1\oplus V_2$ with the operators defined in \eqref{Tag1},
suppose that $1 \leq k < n\ \&\  b_p=0 \ \forall\  p<k \ \&\  b_k\neq 0$.  Then 
$\mathcal J_n = 0$ if and only if 
$$
a_{n-k+1} 
= 
 \left( \tfrac {n-k} {k} \right) c_{n-k+1} 
+ \frac
{\sum_{q=2}^{n-k} (-1)^{q(n-q)}b_{n-q+1} \left[\tbinom{n-1}{q-1} a_{q} - \tbinom{n-1}{q-2} c_{q} \right]}
{ (-1)^{k(n-k)+n} \tbinom{n-1}{n-k} b_k }.
$$
\end {lemma}

\begin {proof}
We know that
$$
\mathcal J_n = 0 \Longleftrightarrow 0=
\sum_{p=1}^n (-1)^{p(n-p)}(-1)^n b_p \left[ (-1) \tbinom{n-1}{n-p} a_{n-p+1} 
+ \tbinom{n-1}{p} c_{n-p+1} \right].
$$
Since $b_k\neq 0$ by assumption, $a_1=0$ by the preceding lemma.  Thus when $p=n$, $a_{n-p+1}=0$.
Also, when $p=n$, $\tbinom{n-1}{p}=0$.
Thus we can change the upper limit of the summation to $n-1$.  Furthermore, since $b_p=0\ \forall\ p<k$, the 
lower limit can be changed to $k$. 
Thus $\mathcal J_n = 0$ if and only if 
$$ 
0 = \sum_{p=k}^{n-1} (-1)^{p(n-p)}b_p \left[ (-1) \tbinom{n-1}{n-p} a_{n-p+1} 
+ \tbinom{n-1}{p} c_{n-p+1} \right].
$$
Removing the term with $p=k$ from the summation, we get an equivalent equation:
\begin {multline*}
 (-1)^{k(n-k)}b_k \tbinom{n-1}{n-k} a_{n-k+1} 
=
 (-1)^{k(n-k)}b_k  \tbinom{n-1}{k} c_{n-k+1} \\
+ 
 \sum_{p=k+1}^{n-1} (-1)^{p(n-p)}b_p \left[ (-1) \tbinom{n-1}{n-p} a_{n-p+1} + \tbinom{n-1}{p} c_{n-p+1} \right].
\end {multline*}
Since we supposed that $1 \leq k < n$, $\tbinom{n-1}{n-k}\neq 0$.  Also $b_k\neq 0$ by supposition.
Thus the previous equation is equivalent to the formula,
\begin {multline*}
 a_{n-k+1} 
= 
 \frac { \tbinom{n-1}{k} c_{n-k+1} } { \tbinom{n-1}{n-k} }\\
+ \frac
{\sum_{p=k+1}^{n-1} (-1)^{p(n-p)}b_p \left[ (-1) \tbinom{n-1}{n-p} a_{n-p+1} + \tbinom{n-1}{p} c_{n-p+1} \right]}
{ (-1)^{k(n-k)}b_k \tbinom{n-1}{n-k} }.
\end {multline*}
Substituting $q=n-p+1$, we get
\begin {multline*}
 a_{n-k+1} 
= 
 \frac { \tbinom{n-1}{k} c_{n-k+1} } { \tbinom{n-1}{n-k} } \\
+ 
\frac
{\sum_{q=2}^{n-k} (-1)^{(n-q+1)(q-1)}b_{n-q+1} \left[ (-1) \tbinom{n-1}{q-1} a_{q} + \tbinom{n-1}{n-q+1} c_{q} \right]}
{ (-1)^{k(n-k)}b_k \tbinom{n-1}{n-k} }.
\end {multline*}
Simplifying, we get
$$
 a_{n-k+1} = 
 \left( \tfrac {n-k} {k} \right) c_{n-k+1} 
+ \frac
{\sum_{q=2}^{n-k} (-1)^{q(n-q)}b_{n-q+1} \left[\tbinom{n-1}{q-1} a_{q} - \tbinom{n-1}{q-2} c_{q} \right]}
{ (-1)^{k(n-k)+n} \tbinom{n-1}{n-k} b_k }.
$$
\end {proof}


\begin {lemma} \label {lem3.6}
If $b_p=0\ \forall\ p < n$ then $V=V_0\oplus V_1\oplus V_2$, with the skew operators defined in \eqref{Tag1}
is $L_n$.
\end {lemma}

\begin {proof}
Suppose $1\leq m \leq n$. We know that
$$
\mathcal J_m = 0 \Longleftrightarrow 0=
\sum_{p=1}^m (-1)^{p(m-p)}(-1)^m b_p \left[ (-1) \tbinom{m-1}{m-p} a_{m-p+1} 
+ \tbinom{m-1}{p} c_{m-p+1} \right].
$$
If $m < n$, we have $b_p = 0\ \forall\ p \leq m$, and $\mathcal J_m = 0$.
So we only need to show that $\mathcal J_n = 0$.
\begin {multline*}
\mathcal J_n =
\sum_{p=n}^n (-1)^{p(n-p)}(-1)^n b_p \left[ (-1) \tbinom{n-1}{n-p} a_{n-p+1} 
+ \tbinom{n-1}{p} c_{n-p+1} \right]
\\ =
(-1)^n b_n \left[ (-1) \tbinom{n-1}{o} a_{1} 
+ \tbinom{n-1}{n} c_{1} \right] = (-1)^{n+1} (n-1) b_n a_1.
\end {multline*}
If $b_n = 0$, then this is zero.  If $b_n\neq 0$, 
we know from Lemma \eqref{lem3.4} that $a_1 = 0$,
which makes the sum zero in all cases.
Thus $\mathcal J_n = 0 \ \forall\ m\leq n$.  Therefore, $V$ is $L_n$.
\end {proof}


Now the following theorem provides necessary and sufficient conditions which identify all possible $L_n$ structures 
on $V=V_0\oplus V_1\oplus V_2$.  This also leads to a corollary which specifies all possible
$L_{\infty}$ structures on this space.

\newpage

\begin {theorem} \label {thm3.7}
Suppose that $V=V_0\oplus V_1\oplus V_2$ has skew operators defined by
\phantom{**********************************}
$
\left\{\ 
\begin {aligned}
l_n(v\otimes w^{\otimes n-1})&=a_n w          \!\quad \forall\ n\geq 1 \\
l_n(w^{\otimes n})&=b_n x                      \,\quad \forall\ n\geq 1 \\
l_n(v\otimes w^{\otimes n-2}\otimes x)&=c_n x  \,\quad \forall\ n\geq 2 \\
\end {aligned}
\right\}.
$
$V$ is $L_n$ if and only if 
\begin {itemize}
\item
$b_p=0\ \forall\ p < n$ or 
\item
$\ \exists\ k$ such that $1 \leq k < n \ \&\ b_p = 0 \ \forall\ p < k \ \&\ b_k \neq 0\ \&\ \ \forall\ 1 \leq m \leq n-k+1$,
$$
\quad a_m =
 \left( \tfrac {m-1} {k} \right) c_{m} 
+ \frac
{\sum_{p=2}^{m-1} (-1)^{p(m+k)}b_{m+k-p} \left[\tbinom{m+k-2}{p-1} a_{p} - \tbinom{m+k-2}{p-2} c_{p} \right]}
{ (-1)^{m(m+k)+1} \tbinom{m+k-2}{m-1} b_k }.
$$
\end {itemize}
\end {theorem}

\begin {proof}
If $b_p=0\ \forall\ p < n$ then
we know from Lemma \eqref{lem3.6} that $V$ is $L_n$.  
So suppose otherwise.  

Then $\ \exists\ k$ such that $1 \leq k < n \ \&\ b_p = 0 \ \forall\ p < k \ \&\ b_k \neq 0$.
Using Lemma \eqref{lem3.6} again, we know that $V$ is $L_k$.  Thus $\mathcal J_q = 0\ \forall\ q \leq k$.
Thus
$V$ is $L_n$ if and only if
$\mathcal J_q = 0\ \forall\ k < q\leq n$, which by Lemma \eqref {lem3.5} is true if and only if 
$$
a_{q-k+1} 
= 
 \left( \tfrac {q-k} {k} \right) c_{q-k+1} 
+ \frac
{\sum_{p=2}^{q-k} (-1)^{p(q-p)}b_{q-p+1} \left[\tbinom{q-1}{p-1} a_{p} - \tbinom{q-1}{p-2} c_{p} \right]}
{ (-1)^{k(q-k)+q} \tbinom{q-1}{q-k} b_k }
$$
$\forall\ k < q\leq n$.
If we substitute $m=q-k+1$, this formula becomes
$$
a_{m} 
= 
 \left( \tfrac {m-1} {k} \right) c_{m} 
+ \frac
{\sum_{p=2}^{m-1} (-1)^{p(m+k-1-p)}b_{m+k-p} \left[\tbinom{m+k-2}{p-1} a_{p} - \tbinom{m+k-2}{p-2} c_{p} \right]}
{ (-1)^{k(m-1)+m+k-1}\tbinom{m+k-2}{m-1} b_k }.
$$
Note that $k < q\leq n \Longleftrightarrow 1 < q-k+1 \leq n-k+1$.
Thus $V$ is $L_n$ if and only if 
$$
a_m =
 \left( \tfrac {m-1} {k} \right) c_{m} 
+ \frac
{\sum_{p=2}^{m-1} (-1)^{p(m+k)}b_{m+k-p} \left[\tbinom{m+k-2}{p-1} a_{p} - \tbinom{m+k-2}{p-2} c_{p} \right]}
{ (-1)^{m(m+k)+1} \tbinom{m+k-2}{m-1} b_k }
$$   
$\forall\ 1 < m \leq n-k+1$.
\end {proof}


\begin {corollary} \label {cor3.8}
Suppose that $V=V_0\oplus V_1\oplus V_2$, with the skew operators defined in \eqref{Tag1} is $L_n$.
Then $V$ is $L_{n+1}$ if and only if 
\begin {itemize}
\item
$b_p=0\ \forall\ p < n+1$ or 
\item
$\ \exists\ k$ such that $1 \leq k < n+1 \ \&\ b_p = 0 \ \forall\ p < k \ \&\ b_k \neq 0$ and
$$
\qquad\ \ a_{n-k+2} =
 \left( \tfrac {n-k+1} {k} \right) c_{n-k+2} 
+ \frac
{\sum_{p=2}^{n-k+1} (-1)^{p(n+2)}b_{n+2-p} \left[\tbinom{n}{p-1} a_{p} - \tbinom{n}{p-2} c_{p} \right]}
{ (-1)^{n(n-k)+1} \tbinom{n}{k-1} b_k }.
$$
\end {itemize}
\end {corollary}


\begin {corollary} \label {cor3.9}
$V=V_0\oplus V_1\oplus V_2$ with the skew operators defined in \eqref{Tag1} is $L_{\infty}$
if and only if 
\begin {itemize}
\item
$b_p=0\ \forall\ p\in \mathbb N$ or 
\item
$\exists\ k\in \mathbb N$ such that $b_p = 0 \ \forall\ p < k \quad \& 
\quad b_k \neq 0 \quad  \& \quad \forall\ m \in \mathbb N$,
$$
a_m =
 \left( \tfrac {m-1} {k} \right) c_{m} 
+ \frac
{\sum_{p=2}^{m-1} (-1)^{p(m+k)}b_{m+k-p} \left[\tbinom{m+k-2}{p-1} a_{p} - \tbinom{m+k-2}{p-2} c_{p} \right]}
{ (-1)^{m(m+k)+1} \tbinom{m+k-2}{m-1} b_k }.
$$
\end {itemize}
\end {corollary}


The following result gives necessary and sufficient conditions under which an $L_3$ structure
on $V=V_0\oplus V_1\oplus V_2$ is a differential graded Lie algebra.

\begin {theorem}
Suppose that $V=V_0\oplus V_1\oplus V_2$ with the operators defined in \eqref{Tag1}
is an $L_3$ structure.  Then
$V$ is a differential graded Lie algebra if and only if $b_1=0$ or $b_2=0$ or $a_2=0$.
\end {theorem}

\begin {proof}
Recall that a d.g. Lie algebra is an $L_2$ structure in which $l_2\circ l_2 \equiv 0$.
In this space, $l_2\circ l_2 \equiv 0 \Longleftrightarrow l_2\circ l_2(v\otimes w^{\otimes 2})=0$,
by the argument used in Lemma \eqref {lem3.2}.
\begin {multline*}
l_2\circ l_2(v\otimes w^{\otimes 2}) \\
= l_2(l_2(v\otimes w)\otimes w) - (-1)^1 l_2(l_2(v\otimes w)\otimes w) + (-1)^0 l_2(l_2(w\otimes w)\otimes v)\\
= 2 l_2(a_2 w\otimes w) + l_2(b_2 x\otimes v) = 2 a_2 b_2 x - b_2 c_2 x = (2 a_2 - c_2) b_2 x.
\end {multline*}
Thus $V$ is a d.g. Lie algebra if and only if $(2 a_2 - c_2) b_2 = 0$. 
Note that when $b_2 = 0$, this condition makes $V$ a d.g. Lie algebra.

To prove the rest, we need to use the fact that $V$ is $L_3$.  Recall that in this space,
$$
\mathcal J_n = 0 \Longleftrightarrow
0=\sum_{p=1}^n (-1)^{p(n-p)}(-1)^n b_p \left[ (-1) \tbinom{n-1}{n-p} a_{n-p+1} + \tbinom{n-1}{p} c_{n-p+1} \right].
$$
Therefore, $V$ is $L_3$ if and only if the following 3 constraints hold:
\begin {alignat*}{2}
0&=a_1 b_1
&&\qquad\text{(since } \mathcal J_1 = 0\text ).\\
0&=b_1 a_2 - b_1 c_2 - b_2 a_1
&&\qquad\text{(since }\mathcal J_2 = 0\text ).\\
0&=b_1 ( a_3 - 2c_3) + b_2 ( 2 a_2 - c_2) + b_3 a_1
&&\qquad\text{(since }\mathcal J_3 = 0\text ).
\end {alignat*}

Suppose $b_2\neq 0$ and $b_1=0$.
By plugging into the second Jacobi constraint, we get $a_1=0$ (since $b_2\neq 0$).
Plugging $a_1=0$ and $b_1=0$ into the third constraint, we get $b_2(2a_2-c_2)=0$, 
which makes $V$ is a d.g. Lie algebra.

So suppose $b_2\neq 0$ and $b_1\neq 0$ and $a_2=0$.
From the first constraint, $a_1 = 0$.
Plugging $a_1=a_2=0$ into the second Jacobi constraint, we get 
$b_1 c_2 = 0$.  Thus $c_2=0$ (since $b_1\neq 0$).
Since $a_2=0$ and $c_2=0$, $(2a_2-c_2)b_2=0$, and $V$ is a d.g. Lie algebra.

Finally, suppose $b_2\neq 0$ and $b_1\neq 0$ and $a_2\neq 0$.
From the first constraint, $a_1 = 0$.
Plugging $a_1=0$ into the second Jacobi constraint, we get
$0=b_1(a_2 -  c_2)$.  Thus $c_2 = a_2$ (since $b_1 \neq 0$). 
Then  $(2 a_2 - c_2) b_2 = (2 a_2 - a_2) b_2 = a_2 b_2 \neq 0$ (since we assumed that both
$a_2$ and $b_2$ are nonzero).  Thus $V$ is not a d.g. Lie algebra in this case.
\end {proof}

\begin {example}
Let $V=V_0 \oplus V_1 \oplus V_2$ where
\begin {alignat*} {6}
&l_1(v) &= 0 \qquad & \qquad &l_2(v \otimes w) &= w \qquad & \qquad &l_3(v \otimes w \otimes w)& &= -w \\
&l_1(w) &= x \qquad & \qquad &l_2(w \otimes w) &= x \qquad & \qquad &l_3(w \otimes w \otimes w)& &= 0  \\
&\      &\   \qquad & \qquad &l_2(v \otimes x) &= x \qquad & \qquad &l_3(v \otimes w \otimes x)& &= 0  
\end {alignat*}
In the notation of \eqref{Tag1}, 
$a_1 = b_3 = c3 = 0$, $b_1 = a_2 = b_2 = c_2 = 1$, and $a_3 = -1$.
This example is $L_3$ since $a_1$, $a_2$, and $a_3$ satisfy the formula in 
Theorem \eqref{thm3.7}.  However, it is {\bf not} a differential graded Lie algebra since 
\begin {multline*}
l_2\circ l_2(v\otimes w\otimes w) = 
l_2(l_2(v\otimes w)\otimes w) -(-1)^1 l_2(l_2(v\otimes w)\otimes w) + (-1)^0 l_2(w\otimes w)\otimes v)\\
= l_2(w\otimes w) + l_2(w\otimes w) + 0 = 2 x.
\end {multline*}
Suppose that we extend this example by defining $l_4$:
\begin {alignat*} {3}
&l_4(v\otimes w^{\otimes 3})=a_4 w \qquad & \qquad
&l_4(w^{\otimes 4})=b_4 x \qquad & \qquad
&l_4(v\otimes w^{\otimes 2}\otimes x)=c_4 x
\end {alignat*}
Then by Corollary \eqref{cor3.8}, this structure is $L_4$ if and only if 
$a_{4} = 3 c_{4}-1$.
\qed
\end {example}

\begin {example}
Let $V=V_0 \oplus V_1 \oplus V_2$ where
\begin {alignat*} {3}
V_0 &= <v> \qquad & \qquad V_1 &= <w> \qquad & \qquad V_2 &= <x>
\end {alignat*}
Define skew linear operators $l_n:V^{\otimes n}\rightarrow V$ of degree $2-n$ by
\begin {alignat*} {3}
&l_n(v\otimes w^{\otimes n-1})=0 \qquad & \quad
&l_n(w^{\otimes n})= x \qquad & \quad
&l_n(v\otimes w^{\otimes n-2}\otimes x)=0 \qquad \ \ \forall\ n \in \mathbb N
\end {alignat*}
In the notation of \eqref{Tag1},
$a_n = c_n = 0$ and $b_n = 1 \ \forall\ n \in \mathbb N$, and it is easy to see that
the conditions of corollary \eqref{cor3.9} are satisfied.  
Therefore, this is an $L_{\infty}$ structure on $V$.
Furthermore, it is a differential graded Lie algebra since $a_2 = 0$.
Note that in this simple example, we have a d.g. Lie algebra in which the 
higher operators ($l_n$ where $n>2$) are also nonzero.
\qed
\end {example}

%
%

\section {$\mathbf{V=V_{-1}\oplus V_0\oplus V_1}$ where each is 1-dimensional}

Now we will investigate a vector space in which the grading is shifted down by one
from the previous space.  
This is an interesting example, since the case $V_{-1}\neq 0$ is
important in rational homotopy theory.
Here, we will see that there is a much more limited selection of possible $L_m$
or $L_{\infty}$ structures, and every $L_3$ structure is a differential graded Lie algebra.

It will be shown that every $L_{\infty}$ structure on $V=V_{-1}\oplus V_0\oplus V_1$ 
(using the skew linear operators listed below)
has one of the following forms:
\begin {itemize}
\item
$b_n = 0\ \forall\ n$ but the constants $a_n$, $c_n$ are arbitrary
\item
$b_m \neq 0\ \text { for some } m$ but all constants $a_n$, $c_n$ are zero
\end {itemize}
Now consider $V=\oplus V_i$ where
\begin {alignat*} {4}
V_{-1} &= <u> \qquad & \quad V_0 &= <v> \qquad & \quad V_1 &= <w> \qquad & \quad V_k &= {0},\ k \notin \{ -1,0,1 \}
\end {alignat*}
Define skew linear operators $l_n:V^{\otimes n}\rightarrow V$ of degree $2-n$ by
\begin {equation} \label {Tag2}
\left\{\ 
\begin {aligned}
l_n(u\otimes v\otimes w^{\otimes n-2}) &= a_n u  \!\qquad \forall\ n\geq 2 \\
l_n(u\otimes          w^{\otimes n-1}) &= b_n v  \qquad \forall\ n\geq 1 \\
l_n(         v\otimes w^{\otimes n-1}) &= c_n w  \!\qquad \forall\ n\geq 1 \\ 
l_n(                  w^{\otimes n  }) &= 0         \quad\qquad \forall\ n \\ 
l_n(u^{\otimes i}\otimes v^{\otimes j}\otimes w^{\otimes n-i-j})&=0 \qquad\quad \forall\ i, j>1
\end {aligned}
\right.
\end {equation}
where $a_n$, $b_n$, $c_n$ are constants.  
Since $a_n$ is defined above only for $n\geq 2$, we can define $a_1=0$ for the sake of convenience.

\begin {remark}The above list includes all possible skew linear operators of degree $2-n$, by the 
following argument:
\begin {itemize}
\item 
$l_n(u^{\otimes i}\otimes v^{\otimes j}\otimes w^{\otimes n-i-j})=0 \quad \forall\ j>1$
because a repeated component of even degree makes a skew operator zero.
\item
$l_n(u^{\otimes i}\otimes v^{\otimes j}\otimes w^{\otimes n-i-j})\subset V_{-i+n-i-j+2-n}
=V_{2-2i-j}$.  This can only be nonzero if $2-2i-j \in \{ -1, 0, 1 \}$.
If $j = 0$, then $2-2i \in \{ -1, 0, 1 \}$ which forces $i=1$.
If $j = 1$, then $2-2i-1 \in \{ -1, 0, 1 \}$ which implies that $i = 0$ or $i = 1$.
\end {itemize}
The action of $l_k$ on the remaining generators is determined by its skewness.
\end {remark}

In order to determine when these operators define an $L_m$ structure on $V$, we need to know
when $\mathcal J_n = 0\ \forall\ n \leq m$.  The following three lemmas establish when $\mathcal J_n = 0$.


\begin {lemma}
Given $V=V_{-1}\oplus V_0\oplus V_1$ with the skew operators defined in \eqref{Tag2},
$\mathcal J_n = 0$ if and only if 
$$
\mathcal J_n(u^{\otimes 2}\otimes w^{\otimes n-2})= 0 \ \ \&\ \
\mathcal J_n(u\otimes v \otimes w^{\otimes n-2}) = 0 \ \ \&\ \
\mathcal J_n(u\otimes w^{\otimes n-1})=0.
$$
\end {lemma}

\begin {proof}
Since $deg(\mathcal J_n) = 3-n$, 
$$
\mathcal J_n(V_{-1}^{\otimes i}\otimes V_0^{\otimes j}\otimes V_1^{\otimes n-i-j}) 
\subset V_{-i+n-i-j + 3 - n} = V_{-2i-j+3}.
$$
Thus
$\mathcal J_n(V_{-1}^{\otimes i}\otimes V_0^{\otimes j}\otimes V_1^{\otimes n-i-j})$ 
can only be nonzero if $-2i-j+3\in \{ -1,0,1 \}$ 
(since these are the nonzero components of $V$).
Equivalently,
$\mathcal J_n(V_{-1}^{\otimes i}\otimes V_0^{\otimes j}\otimes V_1^{\otimes n-i-j})$ 
can only be nonzero if $2i+j\in \{ 2,3,4 \}$. 
But 
$\mathcal J_n(V_{-1}^{\otimes i}\otimes V_0^{\otimes j}\otimes V_2^{\otimes n-i-j})=0$ 
whenever $j>1$ (since $\mathcal J_n$ is skew and a parameter of even degree would be repeated).
Thus
$\mathcal J_n(V_{-1}^{\otimes i}\otimes V_0^{\otimes j}\otimes V_2^{\otimes n-i-j})$ 
can only be nonzero when
($j=0$ and $i=1$) or ($j=0$ and $i=2$) or ($j=1$ and $i=1$). 
\end {proof}


\begin {lemma}
Given $V=V_{-1}\oplus V_0\oplus V_1$ with the skew operators defined in \eqref{Tag2},
$$
\mathcal J_n(u\otimes w^{\otimes n-1}) =\sum_{p=1}^{n} (-1)^{p(n-p)} \binom{n-1}{p-1} b_p c_{n-p+1}w.
$$
\end {lemma}

\begin {proof}
Similar to Lemma \eqref {lem3.3}.
\end {proof}


\begin {lemma}
Given $V=V_{-1}\oplus V_0\oplus V_1$ with the skew operators defined in \eqref{Tag2},
$$
\mathcal J_n(u^{\otimes 2}\otimes w^{\otimes n-2}) =\sum_{p=1}^{n-1} (-1)^{p(n-p)}(-2) \binom{n-2}{p-1} b_p a_{n-p+1}u.
$$
\end {lemma}

\begin {proof}
Before calculating $l_{n-p+1}\circ l_p(u^{\otimes 2}\otimes w^{\otimes n-2})$, 
it is helpful to consider the ways in which a $(p, n-p)$ unshuffle of 
$u^{\otimes 2}\otimes w^{\otimes n-2}$ could possibly rearrange the two $u$ terms.
\begin {itemize}
\item
If both $u$ are on the left,  
$l_{n-p+1}(l_p(u^{\otimes 2}\otimes w^{\otimes p-2})\otimes w^{\otimes n-p})=l_{n-p+1}(0)=0$.
\item
If both $u$ are on the right, 
$l_{n-p+1}(l_p(w^{\otimes p})\otimes u^{\otimes 2}\otimes w^{\otimes n-p-2})=l_{n-p+1}(0)=0$.
\end {itemize}
Thus we only need to consider unshuffles $\sigma$ which put one $u$ on each side.
Because of ordering, these are the unshuffles such that either
($\sigma (1) = 1 \ \&\ \sigma (p+1) = 2$) or 
($\sigma (1) = 2 \ \&\ \sigma (p+1) = 1$). 
\vspace {2pt}

There are $\tbinom{n-2}{p-1}$ unshuffles such that ($\sigma (1) = 1 \ \&\ \sigma (p+1) = 2$).
In each case, $\chi(\sigma)=+1$,
since each of these permutations can be accomplished solely by commuting elements of odd degree.
Similarly, there are $\tbinom{n-2}{p-1}$ unshuffles such that ($\sigma (1) = 2 \ \&\ \sigma (p+1) = 1$),
each with $\chi(\sigma)=+1$.  Thus

\begin {align*}
l_{n-p+1}\circ l_p(u^{\otimes 2}\otimes w^{\otimes n-2}) 
&=
2 \tbinom{n-2}{p-1} (+1) l_{n-p+1}(l_p(u^\otimes w^{\otimes p-1})\otimes u\otimes w^{\otimes n-p+1})
\\
&=
2 \tbinom{n-2}{p-1} l_{n-p+1}(b_p v\otimes u\otimes w^{\otimes n-p+1})
\\
&=
2 \tbinom{n-2}{p-1} b_p (-1) l_{n-p+1}(u\otimes v\otimes w^{\otimes n-p+1})
\\
&=
(-2) \tbinom{n-2}{p-1} b_p a_{n-p+1} u.
\end {align*}
Thus
$$
\mathcal J_n(u^{\otimes 2}\otimes w^{\otimes n-2}) = 
\sum_{p=1}^n (-1)^{p(n-p)}(-2) \tbinom{n-2}{p-1} b_p a_{n-p+1}u. 
$$
Since $a_1=0$, we get 
$$
\mathcal J_n(u^{\otimes 2}\otimes w^{\otimes n-2}) =  \sum_{p=1}^{n-1} (-1)^{p(n-p)}(-2) \tbinom{n-2}{p-1} b_p a_{n-p+1}u. 
$$
\end {proof}


\begin {lemma}
Given $V=V_{-1}\oplus V_0\oplus V_1$ with the skew operators defined in \eqref{Tag2},
$$
\mathcal J_n(u\otimes v \otimes w^{\otimes n-2})
=\sum_{p=1}^{n} (-1)^{p(n-p)}b_{n-p+1} \left[ \binom{n-2}{p-2} a_p - \binom{n-2}{p-1} c_p \right]v.
$$
\end {lemma}

\begin {proof}
Denote $x_1\otimes\dots\otimes x_n = u\otimes v \otimes w^{\otimes n-2}$.
Consider the different ways in which a $(p,n-p)$ unshuffle $\sigma$ could possibly arrange the $u$ and $v$ terms:
\begin {itemize}
\item
If both $u$ and $v$ are on the left,
$l_{n-p+1}(l_p(u\otimes v\otimes w^{\otimes p-2})\otimes w^{\otimes n-p})=
l_{n-p+1}(a_p u\otimes w^{\otimes n-p}) = a_p b_{n-p+1} v$.
There are $\tbinom {n-2}{p-2}$ such unshuffles, with $\chi (\sigma) = +1$ in each case.
\item
If $u$ is on the left and $v$ is on the right,
$l_{n-p+1}(l_p(u\otimes w^{\otimes p-1})\otimes v\otimes w^{\otimes n-p-1})=
l_{n-p+1}(b_p v^{\otimes 2}\otimes w^{\otimes n-p-1})=0$
\item
If $v$ is on the left and $u$ is on the right, 
$l_{n-p+1}(l_p(v\otimes w^{\otimes p-1})\otimes u\otimes w^{\otimes n-p+1})
=l_{n-p+1}(c_p w\otimes u\otimes w^{\otimes n-p+1})
= c_p b_{n-p+1} v$.  There are $\tbinom {n-2}{p-1}$ of these, with $\chi(\sigma )=-1$.
\item
If both are on the right, $l_p(w^p)=0$.
\end {itemize}
Thus
$$
\mathcal J_n(u\otimes v \otimes w^{\otimes n-2})
= \sum_{p=1}^n (-1)^{p(n-p)}  \left[ \tbinom{n-2}{p-2} a_p b_{n-p+1} v - \tbinom{n-2}{p-1} c_p b_{n-p+1} v \right].
$$
\end {proof}
Thus $\mathcal J_n=0$ if and only if the following 3 constraints hold:
\begin {alignat*}{2}
0&=\sum_{p=1}^{n} (-1)^{p(n-p)} \binom{n-1}{p-1} b_p c_{n-p+1} 
&&\qquad\text{(from }u\otimes w^{\otimes n-1}\text ).\\
0&=\sum_{p=1}^{n-1} (-1)^{p(n-p)}\binom{n-2}{p-1} b_p a_{n-p+1}
&&\qquad\text{(from }u^{\otimes 2}\otimes w^{\otimes n-2}\text ).\\
0&=\sum_{p=1}^{n} (-1)^{p(n-p)}b_{n-p+1} \left[ \binom{n-2}{p-2} a_p - \binom{n-2}{p-1} c_p \right]
&&\qquad\text{(from } u\otimes v\otimes w^{\otimes n-2}\text ).
\end {alignat*}


The following lemma shows that if any constant $b_k$ is nonzero, then this forces 
many other constants $a_q, c_q$ to be zero in any $L_m$ structure on $V=V_{-1}\oplus V_0\oplus V_1$.

\begin {lemma}
Suppose $1\leq k\leq m$ and $b_k \neq 0$ and $b_p = 0 \ \forall\ p<k$.
Then
$V=V_{-1}\oplus V_0\oplus V_1$ with the skew linear operators defined in \eqref{Tag2}
is $L_m$ if and only if $a_q = c_q = 0\ \forall\ q\leq m-k+1$.
\end {lemma}

\begin {proof}
Since $b_p = 0\ \forall\ p < k$, $\mathcal J_n=0$ if and only if
\begin {multline*}
0=\sum_{p=k}^{n} (-1)^{p(n-p)} \tbinom{n-1}{p-1} b_p c_{n-p+1} \ \&\ 
0=\sum_{p=k}^{n-1} (-1)^{p(n-p)}\tbinom{n-2}{p-1} b_p a_{n-p+1} \\ \&\ 
0=\sum_{p=1}^{n-k+1} (-1)^{p(n-p)}b_{n-p+1} \left[ \tbinom{n-2}{p-2} a_p - \tbinom{n-2}{p-1} c_p \right].
\end {multline*}
Substituting $q=n-p+1$, this is equivalent to 
\begin {multline*}
\!\!0=\!\!\!\!\sum_{q=1}^{n-k+1}\! (-1)^{(n+1)(q-1)} \tbinom{n-1}{n-q} b_{n-q+1} c_{q}  \quad \& \quad
0=\!\!\!\!\sum_{q=2}^{n-k+1}\! (-1)^{(n+1)(q-1)}\tbinom{n-2}{n-q} b_{n-q+1} a_{q} \\ \&\ 
0=\!\!\sum_{p=1}^{n-k+1} (-1)^{p(n-p)}b_{n-p+1} \left[ \tbinom{n-2}{p-2} a_p - \tbinom{n-2}{p-1} c_p \right].
\end {multline*}

Suppose $a_q = c_q = 0\ \forall\ q\leq m-k+1$ and $n\leq m$.  Then
$\mathcal J_n = 0\ \forall\ n\leq m$ (since all of the above sums are clearly zero in that case).  Thus $V$ is $L_m$.

Now suppose that $V$ is $L_m$.  
Since we assumed $1\leq k\leq m$, $V$ is $\mathcal J_k$.  When we plug $n=k$ into
$$
0=\sum_{q=1}^{n-k+1} (-1)^{(n+1)(q-1)} \tbinom{n-1}{n-q} b_{n-q+1} c_{q}
$$
the expression simplifies to $0 =  b_{k} c_{1}$.
Since $b_k \neq 0$ by assumption, $c_1 = 0$.
Also, we know $a_1=0$ because of the way that the skew operators were defined.

Now we can do induction:  Suppose $2 \leq i \leq m-k+1 \ \&\ a_q = c_q = 0\ \forall\ q < i$.
We will now show that $a_i = c_i = 0$, which will complete the proof!
Since $i\leq m-k+1$, $i+k-1\leq m$.  Thus $\mathcal J_{i+k-1} = 0$.  Thus we can plug $n=i+k-1$ into 
$$
0=\!\!\!\sum_{q=1}^{n-k+1} (-1)^{(n+1)(q-1)} \tbinom{n-1}{n-q} b_{n-q+1} c_{q}  \quad \& \quad\!
0=\!\!\!\sum_{q=2}^{n-k+1} (-1)^{(n+1)(q-1)}\tbinom{n-2}{n-q} b_{n-q+1} a_{q}. 
$$
Since we assumed that $a_q = c_q = 0\ \forall\ q < i$, this simplifies to 
$$
0=(-1)^{(k)(i-1)} \tbinom{i+k-2}{k-1} b_{k} c_{i} \quad \&\ \quad
0=(-1)^{(k)(i-1)}\tbinom{i+k-3}{k-1} b_{k} a_{i}. 
$$
But the binomial expressions are nonzero since $i\geq 2$, and $b_k\neq 0$ by assumption.
Thus $c_i=0$ and $a_i=0$.
\end {proof}


Now we can get to the main results:  Theorems \eqref{thm4.6}  and \eqref{thm4.7} give
necessary and sufficient conditions under which skew linear operators form an $L_m$ structure on 
$V = V_{-1}\oplus V_0 \oplus V_1$.

\begin {theorem} \label {thm4.6}
Suppose
$V=V_{-1}\oplus V_0\oplus V_1$ has skew linear operators defined by
\phantom{********************************}
$
\left\{\ 
\begin {aligned}
l_n(u\otimes v\otimes w^{\otimes n-2}) &= a_n u  \!\qquad \forall\ n\geq 2 \\
l_n(u\otimes          w^{\otimes n-1}) &= b_n v  \qquad \forall\ n\geq 1 \\
l_n(         v\otimes w^{\otimes n-1}) &= c_n w  \!\qquad \forall\ n\geq 1 \\ 
\end {aligned}
\right\}.
$
$V$ is $L_m$ if and only if
\begin {itemize}
\item
$b_p=0\ \forall\ p \leq m$ or 
\item
$\exists\ k\in \mathbb N$\ \ such that \ \ $1\leq k\leq m$ \ \ and \ \ $b_p = 0 \ \forall\ p<k$ 
\ \ and \ \ $b_k\neq 0$ \ \ and \ \ 
$a_q = c_q = 0\ \forall\ q\leq m-k+1$.
\end {itemize}
\end {theorem}

\begin {proof}
$\mathcal J_n=0$ if and only if
\begin {multline*}
0=\sum_{p=1}^{n} (-1)^{p(n-p)} \tbinom{n-1}{p-1} b_p c_{n-p+1} \ \&\ 
0=\sum_{p=1}^{n-1} (-1)^{p(n-p)}\tbinom{n-2}{p-1} b_p a_{n-p+1} \\ \&\ 
0=\sum_{p=1}^{n} (-1)^{p(n-p)}b_{n-p+1} \left[ \tbinom{n-2}{p-2} a_p - \tbinom{n-2}{p-1} c_p \right].
\end {multline*}
If $b_p=0\ \forall\ p\leq m$, then 
$\mathcal J_n = 0 \ \forall\ p\leq m$ (since all of the above sums are clearly zero).
Thus $V$ is $L_m$.

Now suppose that $V$ is $L_m$ \& \ $\exists\ b_{\,\overline{k}\,}\neq 0$ such that $\,\overline{k}\ \leq m$.
Then \ $\exists\ k$ such that $1\leq k\leq m$ \ \&\ 
$b_p = 0 \ \forall\ p<k\ \&\ b_k\neq 0$. 
Thus by the preceding lemma,  
$V$ is $L_m$ if and only if $a_q = c_q = 0\ \forall\ q\leq m-k+1$.
\end {proof}


\begin {theorem} \label {thm4.7}
Suppose $V=V_{-1}\oplus V_0\oplus V_1$ is an $L_{\infty}$ structure with the skew linear operators defined in \eqref{Tag2}.
If $\ \exists\, m\in \mathbb N \text{ such that } b_m\neq 0 \text{ then } a_n = c_n = 0\ \forall\ n$.
\end {theorem}

\begin {proof}
Suppose $b_p=0\ \forall\ p<m\ \&\ b_m \neq 0$.
Then the constraint from checking $u^{\otimes 2}\otimes w^{\otimes n-2}$ becomes
$$
0=\sum_{p=m}^{n-1} (-1)^{p(n-p)}\binom{n-2}{p-1} b_p a_{n-p+1}.
$$
Substituting $q=n-p+1$, we get
$$
0=\sum_{q=2}^{n-m+1}(-1)^{(n-q+1)(q-1)}\binom {n-2}{n-q} b_{n-q+1} a_q.
$$
We have $a_1=0$.  Now do induction... Suppose $a_q = 0\ \forall\ q<k$ where $k\geq 2$.  Then
$$
0=\sum_{q=k}^{n-m+1}(-1)^{(n-q+1)(q-1)}\binom {n-2}{n-q} b_{n-q+1} a_q.
$$
In the particular case when $n=m+k-1$, this becomes
\begin {multline*}
0=\sum_{q=k}^{k}(-1)^{(m+k-q)(q-1)}\binom {m+k-3}{m+k-1-q} b_{m+k-q} a_q\\
=(-1)^{(m)(k-1)}\binom {m+k-3}{m-1} b_{m} a_k.
\end {multline*}
Note that since $k\geq 2$, $m+k-3\geq m-1$.  Also $m-1\geq 0$ since $m\in \mathbb N$.
\vspace {2pt}
Thus $\binom {m+k-3}{m+k-1-q} \neq 0$.  Since $b_m\neq 0$ by assumption, $a_k=0$.
Thus by induction, $a_n=0\ \forall\ n$.

A very similar induction argument using the constraint from $u^{\otimes 2}\otimes w^{\otimes n-2}$
proves that $c_n=0\ \forall\ n$.
\end {proof}

\begin {corollary}
Given the skew linear operators defined in \eqref{Tag2}, 
every $L_{\infty}$ structure on $V=V_{-1}\oplus V_0\oplus V_1$ has one of the following forms:
\begin {itemize}
\item
$b_n = 0\ \forall\ n$ but the constants $a_n$, $c_n$ are arbitrary
\item
$b_m \neq 0\ \text { for some } m$ but all constants $a_n$, $c_n$ are zero
\end {itemize}
\end {corollary}

The following result shows that every $L_3$ structure on $V=V_{-1}\oplus V_0\oplus V_1$ is a 
differential graded Lie algebra.

\begin {theorem}
Suppose $V=V_{-1}\oplus V_0\oplus V_1$ is an $L_3$ structure with the skew linear operators defined in \eqref{Tag2}.
Then this structure is a differential graded Lie algebra.
\end {theorem}

\begin {proof}
Recall that a d.g. Lie algebra is an $L_2$ structure in which $l_2\circ l_2 \equiv 0$.
In this space, 
$l_2\circ l_2=0 \Longleftrightarrow 
l_2\circ l_2(u^{\otimes 2}\otimes w) =
l_2\circ l_2(u\otimes v \otimes w) =
l_2\circ l_2(u\otimes w^{\otimes 2}) =0$.

\begin {multline*}
l_2\circ l_2(u^{\otimes 2}\otimes w) \\
= l_2(l_2(u\otimes u)\otimes w) - (-1)^1 l_2(l_2(u\otimes w)\otimes u)\ + (-1)^0 l_2(l_2(u\otimes w)\otimes u)\\
= l_2(0) + l_2(b_2 v\otimes u) + l_2(b_2 v\otimes u) 
= 0 - 2 l_2(b_2 u\otimes v) = -2 b_2 a_2 u.
\end {multline*}
\begin {multline*}
l_2\circ l_2(u\otimes v \otimes w) \\
= l_2(l_2(u\otimes v)\otimes w) - (-1)^0 l_2(l_2(u\otimes w)\otimes v) + (-1)^1 l_2(l_2(v\otimes w)\otimes u)\\
= l_2(a_2 u\otimes w) - l_2(b_2 v\otimes v) - l_2(c_2 w\otimes u)
= a_2 b_2 v - 0 + c_2 b_2 v.
\end {multline*}
\begin {multline*}
l_2\circ l_2(u\otimes w^{\otimes 2}) \\
= l_2(l_2(u\otimes w)\otimes w) - (-1)^1 l_2(l_2(u\otimes w)\otimes w) + (-1)^0 l_2(l_2(w\otimes w)\otimes u)\\
= 2 l_2(b_2 v\otimes w) + 0 = 2 b_2 c_2 w.
\end {multline*}
Thus $V$ is a d.g. Lie algebra if and only if $b_2 a_2 = 0$ and $b_2 c_2 = 0$.
If $b_2=0$, then clearly $V$ is d.g. Lie.

If $b_2 \neq 0$,
then $\exists\ k$ such that $k\leq 2$ and $b_k \neq 0$ and $b_p = 0 \ \forall\ p < k$.
Because of this and the fact that $V$ is $L_3$,
Theorem \eqref{thm4.6} requires that $a_q = c_q = 0\ \forall\ q \leq 3-k+1$. 
Since $k\leq 2$, this forces $a_2 = c_2 = 0$, which again makes $V$ d.g. Lie.
\end {proof}

%
%

\section {Other $\mathbf {V = V_m \oplus V_{m+1} \oplus V_{m+2}}$ where each is 1-dimensional}

Most of these other cases are rather trivial: 
\begin {itemize}
\item
When $m=-2$ or $m=-3$, we have a degenerate case in which the
only possible nonzero operators are $l_1$ and $l_2$.
Therefore, any $L_2$ structure on $V=V_m\oplus V_{m+1}\oplus V_{m+2}$ is $L_2$ is automatically $L_{\infty}$.
\item
When $m\leq -4$ or $m\geq 2$, this forces $l_n=0\ \forall\ n>1$,
which gives a trivial $L_{\infty}$ structure to any such differential graded vector space.
\end {itemize}
The results for $V=V_1\oplus V_2\oplus V_3$ (the only remaining interesting case)
will be briefly summarized below.  The proofs are all very similar to those in the preceding sections.
Consider $V=\oplus V_i$ where 
\begin {alignat*} {4}
V_1 &= <w> \qquad & \qquad V_2 &= <x> \qquad & \qquad V_3 &= <y> \qquad & \qquad V_k &= {0},\ k \notin \{ 1,2,3 \}
\end {alignat*}
Define skew linear operators $l_n:V^{\otimes n}\rightarrow V$ of degree $2-n$ by
\begin {equation}
\left\{\ 
\begin {aligned}
l_n(w^{\otimes n})&=a_n x                               \!\quad \forall\ n\geq 1 \\
l_n(w^{\otimes n-1}\otimes x) &= b_n y                  \quad \forall\ n\geq 1 \\
l_n(w^{\otimes n-i-j}\otimes x^{\otimes i}\otimes y^{\otimes j})&=0 \!\qquad \ \forall\ i>1,\ j>0
\end {aligned}
\right.
\label {Tag3}
\end {equation}
where $a_n$, $b_n$ are constants.

\begin {theorem} \label {thm6.7}
Suppose that $V=V_1\oplus V_2\oplus V_3$ has skew operators defined by
\phantom{***********************}
$
\left\{\ 
\begin {aligned}
l_n(w^{\otimes n})&=a_n x                               \!\quad \forall\ n\geq 1 \\
l_n(w^{\otimes n-1}\otimes x) &= b_n y                  \quad \forall\ n\geq 1 \\
l_n(w^{\otimes n-i-j}\otimes x^{\otimes i}\otimes y^{\otimes j})&=0 \!\qquad \ \forall\ i>1,\ j>0
\end {aligned}
\right\}.
$
$V$ is $L_n$ if and only if 
\begin {itemize}
\item
$b_p=0\ \forall\ p < n$ or 
\item
$\ \exists\ k$ such that $1 \leq k < n \ \&\ b_p = 0 \ \forall\ p < k \ \&\ b_k \neq 0\ \&$
$$
 \quad a_{m} 
= 
\frac
{- \sum_{q=k+1}^{m+k-2} (-1)^{(m+k-1)(q-k)}\tbinom{m+k-1}{q-1} a_{m} b_q}
{\tbinom{m-k+1}{k-1} b_k}
\ \ \forall\ 1 < m\leq n-k+1.
$$
\end {itemize}

\end {theorem}


\begin {corollary} \label {cor6.9}
$V=V_1\oplus V_2\oplus V_3$ with the skew operators defined in \eqref{Tag3} is $L_{\infty}$
if and only if 
\begin {itemize}
\item
$b_p=0\ \forall\ p\in \mathbb N$ or 
\item
$\exists\ k\in \mathbb N$ such that $b_p = 0 \ \forall\ p < k \quad \& 
\quad b_k \neq 0 \quad  \&$
$$ 
 a_{m} 
= 
\frac
{- \sum_{q=k+1}^{m+k-2} (-1)^{(m+k-1)(q-k)}\tbinom{m+k-1}{q-1} a_{m} b_q}
{\tbinom{m-k+1}{k-1} b_k}\ \forall\ m \in \mathbb N.
$$
\end {itemize}
\end {corollary}

\begin {remark}
Thus $V=V_1\oplus V_2\oplus V_3$ can have just two types of $L_{\infty}$ structure:
\begin {itemize}
\item
$b_p=0\ \forall\ p\in \mathbb N$ but all constants $a_n$ are arbitrarily chosen
\item
If some $b_p$ is nonzero, the remaining constants $b_n$
can be arbitrarily chosen, but all constants $a_n$ are then uniquely determined by those
choices.
\end {itemize}
Furthermore, every $L_3$ structure on 
$V=V_1\oplus V_2\oplus V_3$ is a d.g. Lie algebra:
\end {remark}


\begin {theorem}
Suppose that $V=V_1\oplus V_2\oplus V_3$ with the operators defined in \eqref{Tag3}
is an $L_3$ structure. 
Then $V$ is a differential graded Lie algebra.
\end {theorem}

\end{document}